\documentclass[12pt,reqno]{amsart}
\usepackage{graphicx}
\usepackage{amssymb}
\usepackage{amsthm}
\usepackage{amsfonts}
\usepackage{amsmath}
\usepackage{amstext}

\usepackage{graphicx}
\usepackage{amssymb}
\usepackage{amsthm}
\usepackage{amsfonts}
\usepackage{amsmath}
\usepackage{amstext}

\newtheorem*{thmI}{Theorem\;I}
\newtheorem*{thmII}{Theorem\;II}

\newtheorem*{key}{Key~Lemma}
\newtheorem*{flem}{Fundamental Lemma}
\newtheorem*{att/con}{Attention/Convention}

\newtheorem{thm}{Theorem}[section]
\newtheorem{cor}[thm]{Corollary}

\newtheorem{lem}[thm]{Lemma}

\newtheorem{defi}[thm]{Definition}
\newtheorem{add}{Addendum}
\newtheorem{rem}[thm]{Remark}
\newtheorem{example}[thm]{\textit{Example}}

\newcommand{\B}{\mathbb{B}}
\newcommand{\C}{\mathbb{C}}
\newcommand{\F}{\mathcal{F}}

\newcommand{\M}{\mathbb{M}}

\newcommand{\Q}{\mathbb{Q}}
\newcommand{\R}{\mathbb{R}}
\newcommand{\Z}{\mathbb{Z}}

\newcommand{\s}{\smallskip}

\newcommand{\CO}{\mathcal{O}}
\newcommand{\PPP}{\mathbb{P}}

\newcommand{\biV}{\textit{V}^{\,\C}_*}

\newcommand{\rar}{\rightarrow}
\newcommand{\bS}{\textbf{S}_*}
\newcommand{\bSf}{\textbf{S}_{*/f}}
\newcommand{\bSft}{\textbf{S}_{*/{f_t}}}
\newcommand{\bSl}{\textbf{S}_{*/l}}
\newcommand{\Rp}{\R_*^+}
\newcommand{\Rm}{\R_*^-}
\newcommand{\Rpf}{\R_{*/{f}}^+}
\newcommand{\Rmf}{\R_{*/f}^-}

\newcommand{\Rpml}{\R_{*/l}^{\pm}}
\newcommand{\pa}{\pmb{\alpha}}
\newcommand{\paf}{\pmb{\alpha}_f}

\newcommand{\pg}{\pmb{\gamma}}
\newcommand{\pgf}{\pmb{\gamma}_f}
\newcommand{\pz}{\pmb{\zeta}}
\newcommand{\vv}{\vec{v}}
\newcommand{\vl}{\vec{\lambda}}

\newcommand{\pl}{\pmb{\lambda}}
\newcommand{\plf}{\pmb{\lambda}_f}

\newcommand{\pmu}{\pmb{\mu}}

\newcommand{\pmuf}{\pmb{\mu}_f}

\newcommand{\eps}{\epsilon}

\begin{document}
\title{Equisingularity in $\textbf{\textit{R}}^{\textbf{2}}$ As Morse Stability \\In Infinitesimal Calculus}
\author{Tzee-Char Kuo and Laurentiu Paunescu}
\address{School of Mathematics, University of Sydney, Sydney, NSW, 2006, Australia }
\email{tck@maths.usyd.edu.au; laurent@maths.usyd.edu.au}

\date{}

\maketitle

\begin{abstract}
Two seemingly unrelated problems are intimately connected.

The first is the equsingularity problem in $\R^2$: For an analytic
family $f_t:(\R^2,0)\rar (\R,0)$, when should it be called an
``equisingular deformation"? This amounts to finding a suitable
trivialization condition (as strong as possible) and, of course, a
criterion.

The second is on the Morse stability. We define $\R_*$, which is
$\R$ ``enriched" with a class of infinitesimals. How to generalize
the Morse Stability Theorem to polynomials over $\R_*$?

The space $\R_*$ is much smaller than the space used in
Non-standard Analysis. Our infinitesimals are analytic arcs,
represented by fractional power series, \textit{e.g.},
$x=y^3+\cdots$, $x=y^{5/2}+\cdots$, $x=y^{3/2}+\cdots$, are
infinitesimals at $0\in \R$, in descending orders.

Thus, $p_t(x)\!:=f_t(x,y)\!:=x^4-t^2x^2y^2-y^4$ is a family of
polynomials over $\R_*$. This family is not Morse stable:\;a
triple critical point in $\R_*$ splits into three when $t\not=0$.

In our Theorem\,II, (B) is a trivialization condition which can
serve as a definition for equisingular deformation; (A), and (A')
in Addendum\,\ref{BMorse}, are criteria, using the stability of
``critical points" and the ``complete initial form"; (C) is the
Morse stability (Remark (\ref{MandZ})). Theorem I consists of
weaker conditions (a), (b), (c). The detailed proofs will appear
later.

We were inspired by the intriguing discovery of S. Koike
(\cite{ko}) that the Brian\c{c}on-Speder family, while
blow-analytically trivial, admits no contact order preserving
trivialization. The notion of blow-analytic trivialization must be
modified; (B) and (b) are options.
\end{abstract}

\section{Results.}\label{Results}
As in the Curve Selection Lemma, by a \textit{parameterized arc}
at $0$ in $\R^2$ (resp.\,$\C^2$) we mean a \textit{real} analytic
map germ $\vl: [0, \eps)\rar \R^2$ (resp.\,$\C^2$), $\vl(0)=0$,
$\vl(s)\not\equiv 0$. We call the image set, $\pl\!:=Im(\vl)$, a
(geometric) \textit{arc} at $0$, or the \textit{locus} of $\vl$;
call $\vl$ \textit{a} \textit{parametrization} of $\pl$.

\s

Take $\pl \not=\pmu$. The distance from $P\in \pl$ to $\pmu$ is a
fractional power series in $s\!:=\overline{OP}$,
$dist(P,\pmu)=as^h+\cdots$, where $a>0$, $h\in \Q^+$.

We call $\CO (\pl, \pmu)\!:=h$ the \textit{\textbf{contact order}}
of $\pl$ and $\pmu$. Define $\CO(\pl,\pl)\!:=\infty$.

\s

Let $\bS^1$, or simply $\bS$, denote the set of arcs at $0$ in
$\R^2$. This is called the \textit{enriched} \textit{unit circle}
for the following reason. The tangent half line at $0$, $\pmb{l}$,
of a given $\pl$ can be identified with a point of the unit circle
$\textbf{S}^1$. If $\pmb{\lambda}\not=\pmb{l}$, then
$1<\CO(\pmb{\lambda},\pmb{l})<\infty$. Hence we can regard
$\pmb{\lambda}$ as an ``\textit{infinitesimal}" \textit{at}
$\pmb{l}$, and $\bS$ as $\textbf{S}^1$ \textit{``enriched"} with
infinitesimals.

\s

Let $f:(\R^2,0)\rar (\R,0)$ be analytic. Write $ \biV(f)\!:=\{\pz
\in \bS^3|f(z,w)\equiv 0 \;\text{on}\; \pz\}$, where $\bS^3$
denotes the set of arcs at $0$ in $\C^2(=\R^4)$, and $f(z,w)$ is
the complexification of $f$.

For $\pl\in \bS$, write
$\CO(\pl,\biV(f))\!:=\max\{\CO(\pl,\pz)|\,\pz\in\biV(f)\}$. Define
the \textit{\textbf{f-height}} of $\pl$ by
$h_f(\pl)\!:=\CO(\pl,\biV(f))$. Hence $h_f(\pl)=\infty$ if
$f(x,y)\equiv 0$ along $\pl$.

For $\pl_1$, $\pl_2\in \bS$, define $\pl_1 \thicksim_f \pl_2$
\textit{if and only if} $h_f(\pl_1)=h_f(\pl_2)<\CO(\pl_1,\pl_2)$.
(In fact, $h_f(\pl_1)<\CO(\pl_1,\pl_2)$ implies
$h_f(\pl_1)=h_f(\pl_2)$.) The equivalence class of $\pl$ is
denoted by $\pl_f$.

\s

We call $\pl_f$ an \textit{\textbf{f-truncated arc}}, or simply an
\textit{\textbf{f-arc}}. Write $\bSf:=\bS /\thicksim_f$,
$h(\plf)\!:=h_f(\pl)$.

(Intuitively, once $f$ is given, arcs are ``blurred" so that only
the equivalence classes are ``observable". We were tempted to call
$\plf$ an ``$f$-observable".)

\s

Define the \textit{\textbf{contact} \textbf{order}} of $\pl_f$ and
$\pmu_f$ by: if $\pl_f\not =\pmu_f$,
$\CO(\pl_f,\pmu_f)\!:=\CO(\pl,\pmu)$, $\pl\in \pl_f$, $\pmu\in
\pmu_f$; and $\CO(\pl_f,\pl_f)\!:=\infty$. This is well-defined.
Write $\CO(\pl_f,\biV(f))
\!:=\CO(\pmb{\lambda},\textit{V}^{\C}_*(f))$.

\s

From now on we assume $f(x,y)$ is \textit{\textbf{mini-regular}}
in $x$, that is, regular in $x$ of order $m(f)$, the multiplicity
of $f$. (Thus the positive and negative $x$-directions are not
important.)

\s

Let $\Rp$ (resp.\,$\Rpf$) denote those arcs of $\bS$
(resp.\,$\bSf$) in $y>0$, not tangent to the $x$-axis, and $\Rm$
(resp.\,$\Rmf$) denote those in $y<0$. Write $\R_{*}\!:=\Rp \cup
\Rm$, $\R_{*/f}\!:=\Rpf\cup \Rmf$.

\s

Take $\plf$, $\pmuf\in \Rpf$, or  $\in\Rmf$. Define
$\pmb{\lambda}_f \simeq \pmb{\mu}_f $ (read:\,``bar equivalent")
\textit{if and only if} $\text{either}\; \plf =\pmuf, \;\text{or
else} \; h(\plf)=h(\pmuf)= \CO(\plf,\pmuf)$. Call an equivalence
class an \textit{\textbf{f-bar}}. The one containing $\plf$ is
denoted by $B(\pl_f)$, having \textit{\textbf{height}}
$h(B(\pl_f))\!:=h(\pl_f)$. (See \cite{Kuo-L}, \cite{kuo-par},
\cite{Kur-Pau}.)

If $h(\pl_f)=\infty$ then $B(\pl_f)=\{\pl_f\}$, a singleton, and
conversely.

\s

The given coordinates $(x,y)$ yield a coordinate on each bar of
finite height, as follows.

Take $B$, say in $\R_{*/f}^+$, $h(B)<\infty$. Take
$\pmb{\lambda}\in \pmb{\lambda}_f\in B$ with parametrization
$\vec{\lambda}(s)$. Eliminating $s$ ($s\geq 0$) yields a
\textit{unique} fractional power series (as in \cite{walker})
\begin{equation}\label{representation}
x=\lambda(y)=a_1y^{\frac{n_1}{d}}+a_2y^\frac{n_2}{d}+\cdots, \;
 d\leq n_1<n_2<\cdots, \; (y\geq 0).
\end{equation}
Here all $a_i\in \R$. Let $\lambda_B(y)$ denote $\lambda(y)$ with
all terms $y^e$, $e\geq h(B)$, deleted. Observe that for any
$\pmb{\mu}\in \pmb{\lambda}_f\in B$, $\mu(y)$ has the form
$\mu(y)=\lambda_B(y)+uy^{h(B)}+\cdots$, where $u\in \R$ is
\textit{uniquely} determined by $\pmb{\lambda}_f$. We say
$\pmb{\lambda}_f\in B$ has \textit{\textbf{canonical coordinate}}
$u$, writing $\plf\!:=u$. We call $x=\lambda_B(y)$, which depends
only on $B$, the \textit{\textbf{canonical representation}} of
$B$.

\s

Take $B$, $h(B)<\infty$, and $u=\plf\in B$. Let us write
$$f(\lambda_B(y)+uy^{h(B)}+\cdots, y)\!:=I^B_f(u)y^e+\cdots, \; I^B_f(\plf)\!:=I^B_f(u)\not=0.$$

An important observation is that $e$ depends only on $B$, not on
$\plf$; $I^B_f(u)$ depends only on $\plf$, not on $\pl \in \plf$,
and is a polynomial (Lemma\,(\ref{I}) below). We call
$\i{L}_f(B)\!:=L_f(\plf)\!:=e$ the \textit{\textbf{Lojasiewicz
exponent}} of $f$ on $B$.
\begin{att/con}Not every $u\in \R$ is a canonical coordinate. For example,
$f(x,y)=x^2-y^3$ has a bar $B$ of height $3/2$, and $\pm 1$ are
not canonical coordinates; $I^B_f(u)$ is not a priori defined at
$\pm 1$. Since $I^B_f$ is a polynomial, we shall regard it
\textit{as defined for all} $u\in \R$.

In general, the canonical coordinate identifies $B$ with a copy of
$\R$ minus the real roots of $I^B_f$. Hence $\bar{B}$, the metric
space completion, is a copy of $\R$.

If $ B=\{\plf\}$, a singleton,  we define $I^B_f(\plf)\!:= 0$,
$L_f(\plf)\!:=\infty$.
\end{att/con}

Now, take $l(x,y)\!:=x$, and consider $\bSl$. If
$\nu(y)=ay^e+\cdots$, $a\not=0$, $e\geq 1$, then the $l$-arc
$\pmb{\nu}_l$ can be identified with $(a,e)\in (\R-\{0\})\times
\Q^{+1}$, $\Q^{+1}\!:=\{r\in \Q^+|\,r\geq 1\}$. If $\nu(y)\equiv
0$ then $h(\pmb{\nu}_l)=\infty$; we write $\pmb{\nu}_l
\!:=(0,\infty)$. We call $\mathcal{V}\!:=((\R-\{0\})\times
\Q^{+1})\cup \{(0,\infty)\}(=\Rpml)$ the
\textit{\textbf{infinitesimal value space}}. The given $f$,
mini-regular in $x$, induces a $\mathcal{V}$-valued function
$$f_*: \R_{*/f}\rightarrow \mathcal{V}, \;
f_*(\plf)\!:=(I^B_f(\plf), L_f(\plf))\in \mathcal{V}, \; (\plf\in
B).$$

\s

 Take $z\in \C$. We say $z$ is a $B$-\textit{\textbf{root}} of
$f$ if $f$ has a Newton-Puiseux root of the form
$\alpha(y)=\lambda_B(y)+zy^{h(B)}+\cdots$. The number of such
roots is the \textit{multiplicity} of $z$.

\begin{defi}\label{critical point}
Take $c\!:=\pgf\in B$. If $h(B)<\infty$ and $c\,(\in \R)$ is a
$B$-root of $f_x$, say of multiplicity $k$, we say $\pgf$ is a
(\textit{real}) \textit{\textbf{critical point}} of $f_*$ of
multiplicity $m(\pgf)\!:=k$.

If $B=\{\pgf\}$, and $m(B)\geq 2$, we also call $\pgf$ a critical
point of multiplicity $m(B)-1$.

Call $f_*(c)\!:=f_*(\pgf)\in \mathcal{V}$ the
\textit{\textbf{critical value}} at $\pgf$.

If $f_x$ has complex $B$-root(s), but no real $B$-root, then we
take a \textit{generic} real number $r$, put
$\gamma(y)\!:=\lambda_B(y)+ry^{h(B)}$, and call $\pgf$
\textit{the} real critical point in $B$ with multiplicity
$m(\pgf)\!:=1$. (Convention: For different such $B$, we take
\textit{different} generic $r$.)
\end{defi}
The above is the list of all (real) critical points. (If $f_x$ has
no $B$-root, $B$ yields no critical point.) The number of critical
points is finite (Lemma (\ref{I})).

\s

Now, let $\M$ be the maximal ideal of $\R\{s\}$, furnished with
the point-wise convergence topology, that is, the smallest
topology so that the projection maps
$$\pi_N: \M\longrightarrow \R^N,
\quad a_1s+\cdots +a_Ns^N+\cdots \mapsto (a_1,\cdots,a_N), \quad
N\in \Z^+,$$are continuous. Furnish $\bS$, $\bSf$ with the
quotient topologies by the quotient maps
$$p_*:\M^2-\{0\}\rar \bS,\quad
p_{*/f}:\M^2-\{0\}\rar \bSf.$$

Take $\vl\in \M^2$, and a real-valued function, $\alpha$, defined
near $\vl$. We say $\alpha$ is \textit{analytic} at $\vl$ if
$\alpha =\varphi \circ \pi_N$, $\pi_N$ a projection, $\varphi$ an
analytic function at $\pi_N(\vec{\lambda})$ in $\R^N$. This
defines an analytic structure on $\M^2$. We furnish $\bS$ and
$\bSf$ with the quotient analytic structure.

\s

In the following, let $I$ be a sufficiently small neighborhood of
$0$ in $\R$. We write ``\textit{c}-" for ``continuous",
``\textit{a}-" for ``analytic", ``\textit{c/a}-" for ``continuous
(resp.\,analytic)".

\s

Let $F(x,y;t)$ be a given $t$-parameterized $a$-deformation of
$f(x,y)$. That is to say, $F(x,y;t)$ is real analytic in
$(x,y,t)$, defined for $(x,y)$ near $0\in \R^2$, $t\in I$, with
$F(x,y;0)=f(x,y)$, $F(0,0;t)\equiv 0$. When $t$ is fixed, we also
write $F(x,y;t)$ as $f_t(x,y)$.

\s

In $\bS \times I$ define $(\pl,t)\sim_F(\pl^{\prime},t^{\prime})$
\textit{if and only if} $ t=t^{\prime} \;\text{and}\;
\pl\sim_{f_t}\pl^{\prime}.$ Denote the quotient space by
$\bS\times_F I$. Similarly, $\R_{*}^{\pm}\times_F
I\!:=\R_{*}^{\pm}\times I/\sim_F$.

By a $t$-parameterized \textbf{\textit{c/a}-\textit{deformation}}
of $\plf$ we mean a family of $f_t$-arcs, $\pl_{f_t}$, obtained as
follows. Take a parametrization $\vl(s)$ of $\plf$, and a\,
$c/a$-map: $I\rar \M^2$, $t\mapsto \vl_t$, $\vl_0=\vl$. Then
$\pl_{f_t}\!:=p_{*/{f_t}}(\vl_t)$. This is equivalent to taking
a\, \textit{c/a}-map: $I \rar \bS\times_F I$, $t \mapsto
(\pl_{f_t},t)$. A \textit{\textbf{c/a-deformation}} of a given $B$
is, \textit{by definition}, a family $\{B_t\}$ obtained by taking
any $\plf\in B$, a $c/a$-deformation $\pl_{f_t}$, and then
$B_t\!:=B(\pl_{f_t})$.
\begin{thmI}\label{main} The following three conditions are
equivalent.

(\textbf{a}) Each (real) critical point, $\pgf$, of $f_*$ is
\textit{\textbf{stable along}} $\{f_t\}$ in the sense that $\pgf$
admits a $c$-deformation $\pg_{f_t}$, a critical point of
$(f_t)_*$, such that $m(\pg_{f_t})$, $h(\pg_{f_t})$,
$L_{f_t}(\pg_{f_t})$ are constants. (If $\pgf$ arises from the
generic number $r$, we use the same $r$ for $\pg_{f_t}$.)

(\textbf{b}) There exists a ($t$-level preserving) homeomorphism
$$H: (\R^2\times I, 0\times I)\rightarrow (\R^2\times I, 0\times I),
\quad ((x,y),t)\mapsto (\eta_t(x,y),t),$$ which is bi-analytic off
the $t$-axis $\{0\}\times I$, with the following five properties:

{\quad}(b.1) $f_t(\eta_t(x,y))=f(x,y)$, $t\in I$, (trivialization
of $F(x,y;t)$);

{\quad}(b.2) Given any bar $B$, $\eta_t(\vec{\alpha}(s))$ is
analytic in $(\vec{\alpha},s,t)$, $\vec{\alpha}\in
p_{*/f}^{-1}(B)$ (analyticity on each bar); in particular,
$\eta_t$ is arc-analytic, for any fixed $t$;

{\quad}(b.3) $\CO(\pmb{\alpha},
\pmb{\beta})=\CO(\eta_t(\pmb{\alpha}), \eta_t(\pmb{\beta}))$
(contact order preserving); moreover, $\eta_t(\pmb{\alpha}_f)\in
\bSft$ is well-defined (invariance of truncated arcs).

{\quad}(b.4) The induced mapping $\eta_{t}:B\rightarrow B_t$
extends to an analytic isomorphism:\,$\bar{B}\rightarrow
\bar{B}_t$.

{\quad}(b.5) If $c$ is a critical point of $f_*$, then
$c_t\!=\eta_t(c)$ is one of $(f_t)_*$, $m(c)=m(c_t)$.

(\textit{\textbf{c}}) There exists an \textbf{isomorphism}
$H_*:\R_{*/f}\times I \rightarrow \R_*\times_F I$,
$(\paf,t)\mapsto (\eta_{t}(\paf), t)$, preserving critical points
and multiplicities. That is to say, $H_*$ is a homeomorphism,

{\quad}(c.1) Given $B$, $B_t\!:=\eta_{t}(B)$ is a bar,
$h(B_t)=h(B)$, $m(B_t)=m(B)$;

{\quad}(c.2) The restriction of  $\eta_{t}$ to $B$ extends to an
analytic isomorphism $\bar{\eta}_{t}:\,\bar{B}\rar \bar{B}_t$;

{\quad}(c.2) If $c$ is a critical point of $f_*$, then
$c_t\!:=\eta_{t}(c)$ is one of $(f_t)_*$, $m(c)=m(c_t)$.
\end{thmI}
\begin{thmII}The following three conditions are equivalent.

(\textbf{A}) The function $f_*$ is \textit{\textbf{Morse stable
along}} $\{f_t\}$. That is, every critical point is stable along
$\{f_t\}$, and for critical points $c\in B$, $c^{\,\prime}\in
B^{\prime}$, $f_*(c)=f_*(c^{\,\prime})$ implies
$(f_t)_*(c_t)=(f_t)_*(c_t^{\,\prime})$.

(\textbf{B}) There exists $H$, as in (\textit{\textbf{b}}), with
an additional property:

{\quad}(b.6) If $c$, $c^{\,\prime}$ are critical points,
$f_*(c)=f_*(c^{\,\prime})$, then
$(f_t)_*(c_t)=(f_t)_*(c_t^{\,\prime})$.

(\textbf{C}) There exist an isomorphism $H_*$ as in
(\textit{\textbf{c}}), and an isomorphism $K_*:\mathcal{V}\times
I\rar \mathcal{V}\times I$, such that $K_*\circ (f_*\times
id)=\Phi\circ H_*$, where
$\Phi(\pa_{f_t},t)\!:=((f_t)_*(\pa_{f_t}),t)$.
\end{thmII}
\begin{lem}\label{I}
Let $\{z_1,\cdots,z_q\}$ be the set of $B$-roots of $f$ ($z_i\in
\C$), $h(B)<\infty$. Then
$$I^B_f(x)=a\prod_{i=1}^q(x -z_i)^{m_i},\;
0\not=a\in\R,\,\text{a constant},\,\;m_i\;\text{the multiplicity
of}\;z_i.$$ In particular, $I^B_f(x)$ is a polynomial with real
coefficients.

If $c\!:=\pgf\in B$ is a critical point of $f_*$, then
$\frac{d}{dx}I^B_f(c)=0\not =I^B_f(c)$, and conversely. The
multiplicity of $c$ (as a critical point of the polynomial
$I^B_f(x)$) equals $m(\pgf)$.

The number of critical points of $f_*$ in $\R_{*/f}^{+}$
(resp.\,$\R_{*/f}^-$) is bounded by $m(f)-1$.
\end{lem}
\begin{defi}
The degree of $I^B_f(x)$ is called the
\textit{\textbf{multiplicity}} of $B$, denoted by $m(B)$.

We say $B$ is a \textit{\textbf{polar bar}} if $I^B_f(x)$ has at
least two distinct roots (in $\C$), or $B$ is a singleton with
$m(B)\geq 2$. Call
$\mathcal{I}(f)\!:=\{(B,I^B_f)\,|\,B\;\text{polar}\}$ the
\textit{\textbf{complete initial form}} of $f$.
\end{defi}
\begin{cor}Each critical point belongs to a polar bar; each polar
bar contains at least one critical point.
\end{cor}
We recall Morse Theory. Take an $a$-family of real polynomials
$p_t(x)=a_0(t)x^d+\cdots +a_d(t)$, $a_0(0)\not=0$, $t\in I$, as an
$a$-deformation of $p(x)\!:=p_0(x)$. Let $c_0\in \R$ be a critical
point of $p(x)$, of multiplicity $m(c_0)$. We say $c_0$ is
\textit{stable along} $\{p_t\}$, if it admits a $c$-deformation
$c_t$, $\frac{d}{dx}p_t(c_t)=0$, $m(c_t)=m(c_0)$. (A
$c$-deformation $c_t$, if exists, is necessarily an
$a$-deformation.)
\begin{defi}We say $p(x)$ is \textit{\textbf{Morse and zero stable}} along
$\{p_t\}$ if:

{\quad}(i) Every (real) critical point of $p_0(x)$ is stable along
$\{p_t\}$;

{\quad}(ii) For critical points $c_0$, $c^{\,\prime}_0$,
$p_0(c_0)=p_0(c_0^{\,\prime})$ implies
$p_t(c_t)=p_t(c_t^{\,\prime})$.

{\quad(iii)} If $p_0(c_0)=\frac{d}{dx}p_0(c_0)=0$, then
$p_t(c_t)=\frac{d}{dx}p_t(c_t)=0$.
\end{defi}
\begin{rem}\label{MandZ} Theorem II generalizes a version of the
Morse Stability Theorem: If $p(x)$ is Morse and zero stable along
$\{p_t\}$ then there exist analytic isomorphisms $H, K:\,\R\times
I\rightarrow \R\times I$, such that $K\circ (p\times id)=\Phi\circ
H$, $K(0,t)\equiv 0$, where $\Phi(x,t)\!:=(p_t(x),t)$.

That (a)$\Rightarrow$(c) reduces to the following. Given
$x=f_i(t)$, $1\leq i\leq N$, analytic, $f_i(t)\not =f_j(t)$, for
$i\not=j$, $t\in I$. There exists an analytic isomorphism
$H:\R\times I\rightarrow \R\times I$, $(x,t)\mapsto
(\eta_t(x),t)$, $\eta_t(f_i(t))=const$, $1\leq i\leq N$. (Proved
by Cartan's Theorem A, or Interpolation.)
\end{rem}
We say $\mathcal{I}(f)$ is \textit{\textbf{Morse and zero stable}}
along $\{f_t\}$ if each polar $B$ admits a $c$-deformation $B_t$,
a polar bar of $f_t$, such that two of $h(B_t)$, $m(B_t)$,
$L_{f_t}(B_t)$ are constants (we can then show all three are), and
$\{I^{B}_f\}$ is Morse and zero stable along $\{I^{B_t}_{f_t}\}$,
for each $B$.
\begin{add}\label{BMorse}(\textit{\textbf{B}}) is also equivalent to
($\textit{\textbf{A}}^{\,\prime}$): $\mathcal{I}(f)$ is Morse and
zero stable along $\{f_t\}$.
\end{add}
\begin{example}\label{acondition} For $f_t(x,y)=x^3+3tx^2y+3t^2xy^2+t^3y^3-y^4$,
$f(x,y)$ has critical point $\pgf$, $\gamma(y)\equiv 0$, with
deformation $\gamma_t(y)=-ty$, found by a Tschirnhausen transform,
satisfying (\textit{\textbf{A}}). However, for
$g_t(x,y)=x^3+3tx^2y+t^3y^3-y^4$, terms involving $t$ below the
Newton Polygon of $f$ cannot be cleared, (\textit{\textbf{a}}) is
not satisfied. This idea is elaborated in \S\ref{Newton}.
\end{example}

\section{Relative Newton Polygons.}\label{Newton}

Take $\pl$, say in $\R_*^+$, with $\lambda(y)$ . Let us change
variables: $X\!:=x-\lambda(y), \; Y\!:=y$,
$$\mathcal{F}(X,Y)\!:=f(X+\lambda(Y),Y)
\!:=\sum a_{ij}X^iY^{j/d}, \quad i, j \geq 0,\; i+j> 0.$$

In the first quadrant of a coordinate plane we plot a dot at
$(i,j/d)$ for each $a_{ij}\not=0$, called a (Newton) dot. The
Newton polygon of $\mathcal{F}$ in the usual sense is called the
\textit{Newton Polygon} \textit{of} $f$ \textit{relative to}
$\pl$, denoted by $\PPP (f,\pl)$. (See \cite{kuo-par}.) Write
$m_0\!:=m(f)$. Let the vertices be
$$V_0=(m_0,0),\dots,
V_k=(m_k,q_k),\; q_i\in \Q^+,\, m_i>m_{i+1},\,q_i<q_{i+1}.$$

The (Newton) \textit{edges} are: $E_i=\overline{V_{i-1}V_i}$, with
\textit{angle} $\theta_i$, $\tan
\theta_i\!:=\frac{q_i-q_{i-1}}{m_{i-1}-m_i}$, $\pi/4\leq \theta_i<
\pi/2$; a vertical one, $E_{k+1}$, sitting at $V_k$,
$\theta_{k+1}=\pi/2$; a horizontal one, $E_0$, which is
unimportant.

If $m_k\geq 1$ then $f\equiv 0$ on $\pl$. If $m_k\geq 2$, $f$ is
singular on $\pl$. If $\pl \sim_f \pl^{\prime}$ then
$\PPP(f,\pl)=\PPP(f,\pl^{\prime})$, hence $\PPP(f,\plf)$ is
well-defined.

\textbf{Notation}: $L(E_i)\!:=\overline{V_{i-1}V_i^{\prime}}$,
$V_i^{\prime}\!:=(0,q_{i-1}+m_{i-1}\tan \theta_i)$, \textit{i.e.}
$E_i$ extended to the $y$-axis.

\begin{flem}\label{flem}
Suppose each polar bar $B$ admits a $c$-deformation $B_t$ such
that $h(B_t)$ and $m(B_t)$ are independent of $t$. Then each
$\pl_f\in \R_{*/f}$ admits an \textit{a}-deformation $\pl_{f_t}\in
\R_{*/f_t} $ such that $\PPP(f_t,\pl_{f_t})$ is independent of
$t$. The induced deformation $B_t\!:=B(\pl_{f_t})$ of
$B_0\!:=B(\plf)$, and hence the $a$-deformation
$x=\lambda_{B_t}(y)$ of the canonical representation
$x=\lambda_{B_0}(y)$, are uniquely defined; that is, if we take
any $\pmb{\eta}_f\in B(\plf)$, and a \textit{c}-deformation
$\pmb{\eta}_{f_t}$ with
$\PPP(f_t,\pmb{\eta}_{f_t})=\PPP(f,\pl_f)$, then
$B(\pmb{\eta}_{f_t})=B(\pl_{f_t})$.

Given $B$, $B^{\prime}$. The contact order
 $\CO(B_t,B_t^{\prime})$, defined below, is independent of $t$.
\end{flem}
For $B\not=B^{\prime}$, define
$\CO(B,B^{\prime})\!:=\CO(\pl_f,\pl^{\prime}_f)$, $\pl_f\in B$,
$\pl^{\prime}_f\in B^{\prime}$; and $\CO(B,B)\!:=\infty$.
\begin{example}
For $x^2+2xy-ty^2$, obviously equisingular, the usual Newton
Polygon depends of $t$. This shows the relevance of
\textit{merely} considering Polygons relative to critical points.
\end{example}
The Lemma is proved by a succession of Tschirnhausen transforms at
the vertices, beginning at $V_0$, which represents $a_{m0}X^m$ in
$\mathcal{F}(X,Y)$, $m\!:=m(f)$. Let us define $\mathcal{P}$ by
\begin{equation}\label{perturbation}
F(X+\lambda(Y),Y;t)\!:=\mathcal{F}(X,Y) +\mathcal{P}(X,Y;t),\;
\mathcal{P}(X,Y;t)\!:=\sum p_{ij}(t)X^iY^{j/d}, \end{equation}
where $p_{ij}(t)$ are analytic, $p_{ij}(0)=0$. Take a root of
$\frac{\partial^{m-1}}{\partial
X^{m-1}}[a_{m0}X^m+\mathcal{P}(X,Y;t)]=0$,
$$X=\rho_t(Y)\!:=\sum
b_j(t)Y^{j/d},\, b_j(0)=0,\; b_j(t) \; \text{analytic}.\,
(\text{Implicit Function Theorem})$$ Thus, $\lambda(y)+\rho_t(y)$
is an $a$-deformation of $\lambda(y)$. Let $X_1\!:=X-\rho_t(Y)$,
$Y_1\!:=Y$. Then
$$F(X_1+\lambda(Y_1)+\rho_t(Y_1),Y_1;t)
\!:=\mathcal{F}(X_1,Y_1)+\mathcal{P}^{(1)}(X_1,Y_1;t),$$ where
$\mathcal{P}^{(1)}\!:=\sum p_{ij}^{(1)}(t)X_1^iY_1^{j/d},$
$p_{ij}^{(1)}(0)=0$, and
 $p_{m-1,j}^{(1)}(t)\equiv 0$ (Tschirnhausen).

\s
 For brevity, we shall write the coordinates $(X_1,Y_1,t)$ simply
 as $(X,Y,t)$, abusing notations. That is, we now have $p_{m-1,j}(t)\equiv 0$
  in (\ref{perturbation}).

We claim that $\mathcal{P}$ in fact has no dot below $L(E_1)$.
This is proved by contradiction.

Suppose it has. Take a generic number $s\in \R$. Let
$\zeta(y)\!:=\lambda(y)+sy^e$, $e\!:=\tan \theta_1$, and
$$F(\widetilde{X}+\zeta(\widetilde{Y}), \widetilde{Y};t)\!:
=\mathcal{F}(\widetilde{X},\widetilde{Y})+\widetilde{\mathcal{P}},
\quad \widetilde{\mathcal{P}}(\widetilde{X},\widetilde{Y};0)\equiv
0.$$

Since $s$ is generic, $\PPP(f,\pz_f)$ has only one edge, which is
$L(E_1)$, and $B(\pz_f)$ is polar. Below $L(E_1)$,
$\widetilde{\mathcal{P}}$ has at least one dot (when $t\not=0$),
but still no dot of the form $(m-1,q)$.

A $c$-deformation $B_t$ of $B(\pz_f)$ would either create new
dot(s) of the form $(m-1,q)$ below $L(E_1)$, or else not change
the existing dot(s) of $\widetilde{\mathcal{P}}$ below $L(E_1)$.
(This is the spirit of the Tschirnhausen transformation.) Thus, as
$t\not=0$, $h(B_t)$ or $m(B_t)$, or both, will drop. This
contradicts to the hypothesis of the Fundamental Lemma.

\s

This argument can be repeated recursively at $V_1$, $V_2$, etc.,
to clear all dots under $\PPP(f,\pmb{\lambda}_f)$. More precisely,
suppose in (\ref{perturbation}), $\mathcal{P}$ has no dots below
$L(E_i)$, $0\leq i\leq r$. By the Newton-Puiseux Theorem, there
exists a root $\rho_t$ of $\frac{\partial^{m_r-1}}{\partial
X^{m_r-1}}[aX^{m_r}Y^{q_r}+\mathcal{P}]=0$ with $\CO_y(\rho_t)\geq
\tan\theta_{r+1}$, where $aX^{m_r}Y^{q_r}$ is the term for $V_r$.
A Tschirnhausen transform will then eliminate all dots of
$\mathcal{P}$ of the form $(m_r-1, q)$. As before, all dots below
$L(E_{r+1})$ also disappear.

\s

We have seen the \textit{only} way to clear dots below $\PPP
(f,\pmb{\lambda}_f)$ is by the Tschirnhausen transforms. If
$\PPP(f,\pmb{\eta}_{f_t})=\PPP (f,\pmb{\lambda}_f)$, we must have
$\CO(\pmb{\lambda}_{f_t}, \pmb{\eta}_{f_t})\geq h(B_0)$. The
uniqueness follows.

\s

Define a partial ordering ``$>$" by: $B>\hat{B} $ if and only if $
h(B)>h(\hat{B})=\CO({\pmb{\lambda}_f,\pmb{\mu}_f}),\,
\pmb{\lambda}_f\in B,\pmb{\mu}_f\in \hat{B}$. Let $\hat{B}$ be the
largest bar so that $B\geq \hat{B}$, $B^{\prime}\geq\hat{B}$. We
write $\lambda_B(y)=\lambda_{\hat{B}}(y)+ay^e+\cdots$,
$\lambda_{B^{\prime}}(y)=\lambda_{\hat{B}}(y)+by^e+\cdots$,
$e\!:=h(\hat{B})$. The uniqueness of $\hat{B}_t$ completes the
proof.
\section{Vector fields.}\label{v.f.}

Assume ($\textit{\textbf{a}}$). We use a vector field $\vec{v}$ to
prove $(\textbf{\textit{b}})$. The other implications are not
hard.

Take a critical point $\pgf$, say in $B$,
$\gamma(y)=\lambda_B(y)+cy^{h(B)}$. Let $B_t$ be the deformation
of $B$. Let $c_t$ be the $a$-deformation of $c$,
$\frac{d}{dx}I^{B_t}_{f_t}(c_t)=0$, $m(c_t)=m(c)$. (If $c$ is
generic, take $c_t=c$.)

Let $\gamma_t(y)\!:=\lambda_{B_t}(y)+c_ty^{h(B_t)}$. Then $\pg_t$
is a critical point of $f_t$ in $B_t$.

\s

Now, let $\pgf^{(i)}$, $1\leq i\leq N$, denote all the critical
points of $f$, for \textit{all} (polar) $B$. For brevity, write
$\pg^{(i)}\!:=\pgf^{(i)}$, with deformations $\pg^{(i)}_t$, just
defined.

\s

We can assume $F(x,0;t)=\pm x^m$, and hence $\frac{\partial
F}{\partial t}(x,0;t)\equiv 0$. As $F(x,0;t)=a(t)x^m+\cdots,\,
a(0)\not=0$, a substitution $u=\sqrt[m]{|a(t)|}\cdot x+\cdots$
will bring $F(x,0,t)$ to this form.

We can also assume $\pg^{(i)}\in \Rpf$ for $1\leq i\leq r$, and
$\pg^{(i)}\in\Rmf$ for $r+1\leq i\leq N$.

\s

For each $\pg^{(i)}\in \Rpf$, we now construct a vector field
$\vv_i^+(x,y,t)$, defined for $y\geq 0$.

Write $\pg_t\!:=\pg_t^{(i)}$. Let $X\!:=x-\gamma_t(y), Y\!:=y$.
Then $\mathcal{F}(X,Y;T)\!:=F(X+\gamma_t(Y),Y;T)$ is analytic in
$(X,Y^{1/d},T)$. As in \cite{F-Y}, \cite{Pau}, define
$\vec{v}^{\,+}_i(x,y,t)\!:=\vec{V}(x-\gamma_t(y),y,t)$, $y\geq 0$,
where
\begin{equation}\label{V}
\vec{V}(X,Y,t)\!:=\frac{X\mathcal{F}_X\mathcal{F}_t}
{(X\mathcal{F}_X)^2+(Y\mathcal{F}_Y)^2}\cdot
X\frac{\partial}{\partial X}+\frac{Y\mathcal{F}_Y\,\mathcal{F}_t}
{(X\mathcal{F}_X)^2+(Y\mathcal{F}_Y)^2}\cdot
Y\frac{\partial}{\partial Y}-\frac{\partial}{\partial t}.
\end{equation}

In general, given $\pmb{\alpha}_i$, $x=\alpha_i(y)$, say in $
\R_{*}^+$, $1\leq i \leq r$. Let
$q(x,y)\!:=\prod_{k=1}^r(x-\alpha_k(y))^2$,
$$q_i(x,y)\!:=q(x,y)/(x-\alpha_i(y))^2, \quad p_i(x,y)\!:=q_i(x,y)/[q_1(x,y)+\cdots +q_r(x,y)].$$ We
call $\{p_1,\cdots,p_r\}$ a \textit{\textbf{partition of unity}}
for $\{\pa_1,\cdots,\pa_r\}$.

\s

Now, take $\{p_i\}$ for $\{\pg^{(1)}_t, \cdots \pg^{(r)}_t\}$.
Define $\vec{v}^{\,+}(x,y,t)\!:=\sum_{i=1}^r p_i(x,y,t)\,
\vec{v}_i^{\,+}(x,y,t)$.

Similarly, $\pg^{(i)}_f$, $r+1\leq i \leq N$, yield
$\vec{v}^{\,-}(x,y,t)$, $y\leq 0$. We can then glue
$\vec{v}^{\,\pm}(x,y,t)$ together along the $x$-axis, since
$\vec{v}^{\, \pm}(x,0,t)\equiv -\frac{\partial}{\partial t}$. This
is our vector field $\vec{v}(x,y,t)$, which, by (\ref{V}), is
clearly tangent to the level surfaces of $F(x,y;t)$, proving
($b.1$).

\section{Sketch of Proof}\label{proof}

\begin{lem}\label{Euler}
Let $W(X,Y)$ be a weighted form of degree $d$, $w(X)=h$, $w(Y)=1$.
Take $u_0$, not a multiple root of $W(X,1)$. If $W(u_0,1)\not=0$
or $u_0\not =0$ then, with $X=uv^h$, $Y=v$,
$$|XW_X|+|YW_Y|=unit\cdot \mid\!v\!\mid^d,\;\text{for}\;\, u\;\text{near}\;\, u_0.$$
\end{lem}
For, by Euler's Theorem, if $X-u_0Y^h$ divides $W_X$ and $W_Y$,
then $u_0$ is a multiple root.

To show $(b.2)$, etc., take $\pa$, say in $\Rp$. Take $k$,
$\CO(\pg^{(k)},\pa)=\max\{\CO(\pg^{(j)},\pa)|1\leq j\leq r\}$.

We can assume $\pa$ is not a multiple root of $f$,
$e\!:=\CO(\pg^{(k)},\paf)<\infty$. (If $\pa$ is, then
$\pg^{(k)}=\paf$, $h(B)=\infty$. This case is easy.)

Write $B\!:=B(\paf)$ if $B(\paf)\leq B(\pg^{(k)})$, and
$B\!:=B(\pg^{(k)})$ if $B(\paf)>B(\pg^{(k)})$.

Thus $\alpha(y)=\lambda_B(y)+ay^e+\cdots$,
$\frac{d}{du}I^B_f(a)\not=0$. Let us consider the mapping
$$\tau : (u,v,t)\mapsto
(x,y,t)\!:=(\lambda_{B_t}(v)+uv^e,v,t),\; u\in \R,\;0\leq v
<\varepsilon, \; t\in I,$$ $B_t$ the deformation of $B$, and the
liftings $\vec{\nu}_j^+\!:=(d\tau)^{-1}(p_j\vec{v}_j^+)$,
$\vec{\nu}^+\!:=\sum_{j=1}^r\vec{\nu}_j^+.$
\begin{key}\label{key}
The lifted vector fields $\vec{\nu}_j^+$, and hence $\vec{\nu}^+$,
are analytic at $(u,v,t)$, if $u$ is not a multiple root of
$I^{B_t}_{f_t}$. Moreover, $\vec{\nu}^+(u,0,t)$ is analytic for
all $u\in \R$; that is, $\lim_{v\rightarrow
0^+}\vec{\nu}^+(u,v,t)$ has only removable singularities on the
$u$-axis.
\end{key}
We analyze each $\vec{\nu}_i^+$, using (\ref{V}). For brevity,
write $\B\!:=B(\pg^{(i)})$, $\B_t\!:=B(\pg^{(i)}_t)$.

First, consider the case $B=\B$. This case exposes the main ideas.

Now $I^B_f$ and $\PPP(f,\pg^{(i)})$ are related as follows. Let
$W(X,Y)=\sum_{i,j}a_{ij}X^iY^{j/d}$ be the (unique) weighted form
such that $W(u,1)=I^B_f(u+c)$, $w(X)=h(B)$, $w(Y)=1$, where $c$ is
the canonical coordinate of $\pg^{(i)}$. The Newton dots on the
highest compact edge of $\PPP(f,\pg^{(i)})$ represent the non-zero
terms of $W(X,Y)$; the highest vertex is $(0,L_f(B))$.

Thus $\frac{d}{du}W(0,1)=\frac{d}{du}I^B_f(c)=0$, $W(0,1)\not=0$.
The weighted degree of $W(X,Y)$ is $L_f(B)$.

Hence, by Lemma (\ref{Euler}), the substitution
$X=x-\lambda_B(y)-cy^{h(B)}=(u-c)v^{h(B)}$, $Y=v$, yields
$\CO_v(|X\mathcal{F}_X|+|Y\mathcal{F}_Y|)= L_f(\B)$, if $u-c$ is
not a multiple root of $W(u,1)$.

The Newton Polygon is independent of $t$:\!
$\PPP(f,\pg^{(i)})=\PPP(f_t,\pg_t^{(i)})$. All Newton dots of
$\F$, and hence those of $\mathcal{F}_T$, are contained in
$\PPP(f,\pg^{(i)})$. Hence
$\CO_v(\mathcal{F}_T((u-c)v^{h(B)},v;T))\geq L_f(B)$.

By the Chain Rule, we have $X\frac{\partial}{\partial
X}=(u-c)\frac{\partial}{\partial u}$, $Y\frac{\partial}{\partial
Y}=v\frac{\partial}{\partial v}-h(B)(u-c)\frac{\partial}{\partial
u}$.

It follows that $(d\tau)^{-1}(\vec{v}_i^+)$ and $\vec{\nu}_i$ are
analytic at $(u,v,t)$, if $u$ is not a multiple root of
$I^{B_t}_{f_t}$.

Next, suppose $B<\B$. Again we show $(d\tau)^{-1}(\vec{v}_i^+)$
has the required property.

Write
$\gamma^{(i)}(y)\!:=\lambda_B(y)+c^{\,\prime}y^{h(B)}+\cdots$. Let
$W(X,Y)$ denote the weighted form such that
$W(u,1)=I^B_f(u+c^{\,\prime})$, $w(X)=h(B)$, $w(Y)=1$.

If $W(X,Y)$ has more than one terms, they are dots on a compact
edge of $\PPP(f,\pg^{(i)})$, not the highest one. If $W(X,Y)$ has
only one term, it is a vertex, say $(\bar{m}, \bar{q})$,
$\bar{m}\geq 2$.

In either case, $u=0$ is a multiple root of $W(u,1)$. All Newton
dots of $\mathcal{F}_T$ are contained in $\PPP(f,\pg^{(i)})$. The
rest of the argument is the same as above.

Finally, suppose $B\not \leq\B$. Here $p_i$ plays a vital role in
analyzing $\vec{\nu}_i^+$.

Let $\bar{B}$ denote the largest bar such that $B>\bar{B}\leq \B$.

Let $U\!:=x-\lambda_{B_t}(y)$, $V\!:=y$. The identity
$p_i={p_kq_i}/{q_k}$, and the Chain Rule yield
$$p_i\cdot X\frac{\partial}{\partial X}=p_k\frac{(U+\varepsilon)^2}
{(U+\delta)^2}(U+\delta)\frac{\partial}{\partial U},\quad p_i\cdot
Y\frac{\partial}{\partial
Y}=p_k\cdot\frac{(U+\varepsilon)^2}{(U+\delta)^2}[V\frac{\partial}{\partial
V} -V\delta^{\prime}(V)\frac{\partial}{\partial U}],$$ where
$\delta\!:=\delta(y,t)\!:=\lambda_{B_t}(y)-\gamma_t^{(i)}(y)$,
$\varepsilon\!:=\lambda_{B_t}(y)-\gamma^{(k)}_t(y)$,
$\CO_y(\delta)=h(\bar{B})<h(B)\leq \CO_y(\varepsilon)$.

The substitution $U=uv^{h(B)}$, $V=v$ lifts both to analytic
vector fields in $(u,v,t)$.

It remains to study
$\Psi\!:=\mathcal{F}_T/(|X\mathcal{F}_X|+|Y\mathcal{F}_Y|)$ when
$X=\delta(v,t)+uv^{h(B)}$, $Y=v$.

Let $\mathcal{G}(U,V,T)\!:=\mathcal{F}(U+\delta(V,T),V,T)$. The
Chain Rule yields
\begin{equation}\label{Dominating}
X\mathcal{F}_X=(U+\delta)\mathcal{G}_U,\;
Y\mathcal{F}_Y=V(\mathcal{G}_V-\delta_V \mathcal{G}_U)
,\;\mathcal{F}_T=\mathcal{G}_T-\delta_T\mathcal{G}_U.
\end{equation}
Let us compare $\PPP(f,\pg^{(i)})$ and $\PPP(\mathcal{G}, U=0)$,
the (usual) Newton Polygon of $\mathcal{G}$. Let $E^{\,\prime}_i$,
$\theta_i^{\,\prime}$ and $V_i^{\,\prime}$ denote the edges,
angles and vertices of the latter. Then $E_i=E_i^{\prime}$, for
$1\leq i\leq l$, where $l$ is the largest integer such that $\tan
\theta_l<h(\bar{B})$. Moreover,
\textrm{$\theta_{l+1}^{\,\prime}=\theta_{l+1}$ (although
$E_{l+1}$, } $E^{\,\prime}_{l+1}$ may be different).

Consider the vertex
$V^{\,\prime}_{l+1}\!:=(m_{l+1}^{\prime},q_{l+1}^{\prime})$,
$m_{l+1}^{\prime}\geq 2$. It yields a term $\mu\!:=a(T)U^pV^q$ of
$\delta \mathcal{G}_U$, $a(0)\not=0$, $p\!:=m_{l+1}^{\prime}-1$,
$q\!:=q_{l+1}^{\prime}+\tan \theta_{l+1}$. With the substitution
$U=uv^{h(B)}$, ($u\not=0$,) $V=v$, $\mu$ is the dominating term in
(\ref{Dominating}). That is, $\CO_v(\mu)<\CO_v(\mu^{\prime})$, for
all terms $\mu^{\prime}$ in $U\mathcal{G}_U$, $V\mathcal{G}_V$,
etc., (and for all terms $\mu^{\prime}\not=\mu$ in $\delta
\mathcal{G}_U$), since $\CO_Y(\delta)=\tan\theta_{l+1}$.

It follows that $\Psi$ is analytic. That $\lim \vec{\nu}^+_i$ has
only removable singularities also follows.

Conditions $(b.2)$ etc. can be derived from the Key Lemma.

\end{document}